\documentclass[letterpaper, 10 pt, conference]{ieeeconf}  
     \usepackage{tcolorbox}                                                 
\usepackage{mystyle}
\usepackage{mymacros}

\IEEEoverridecommandlockouts                              
\title{\LARGE \bf
Interpolation Constraints for Computing Worst-Case Bounds in Performance Estimation Problems
}

\author{Anne Rubbens, Nizar Bousselmi, Sébastien Colla, and Julien M. Hendrickx
\thanks{
The authors are with ICTEAM Institute, UCLouvain, Belgium.
{\texttt{ firstname.name@uclouvain.be}}
Their work is supported by FRIA, FNRS, and a L’Oréal-UNESCO For Women in Science fellowship.
This work is also supported
by the ``RevealFlight'' ARC at UCLouvain, by the Incentive Grant
for Scientific Research (MIS) ``Learning from Pairwise Data'' and by the KORNET project from F.R.S.-FNRS
}}

\begin{document}

\maketitle
\thispagestyle{plain}
\pagestyle{plain}
\begin{abstract}
The Performance Estimation Problem (PEP) approach consists in computing worst-case performance bounds on optimization algorithms by solving an optimization problem: one maximizes an error criterion over all initial conditions allowed and all functions in a given class of interest. The maximal value is then a worst-case bound, and the maximizer provides an example reaching that worst case. This approach was introduced for optimization algorithms but could in principle be applied to many other contexts involving worst-case bounds. The key challenge is the representation of infinite-dimensional objects involved in these optimization problems such as functions, and complex or non-convex objects as linear operators and their powers, networks in decentralized optimization etc. This challenge can be resolved by interpolation constraints, which allow representing the effect of these objects on vectors of interest, rather than the whole object, leading to tractable finite dimensional problems. We review several recent interpolation results and their implications in obtaining of worst-case bounds via PEP.
\end{abstract}

\section{Introduction}

\subsection{Worst-case bounds}
\begin{figure}[h]
    \vspace*{-3mm}
    \centering
    \includegraphics[width=0.45\textwidth]{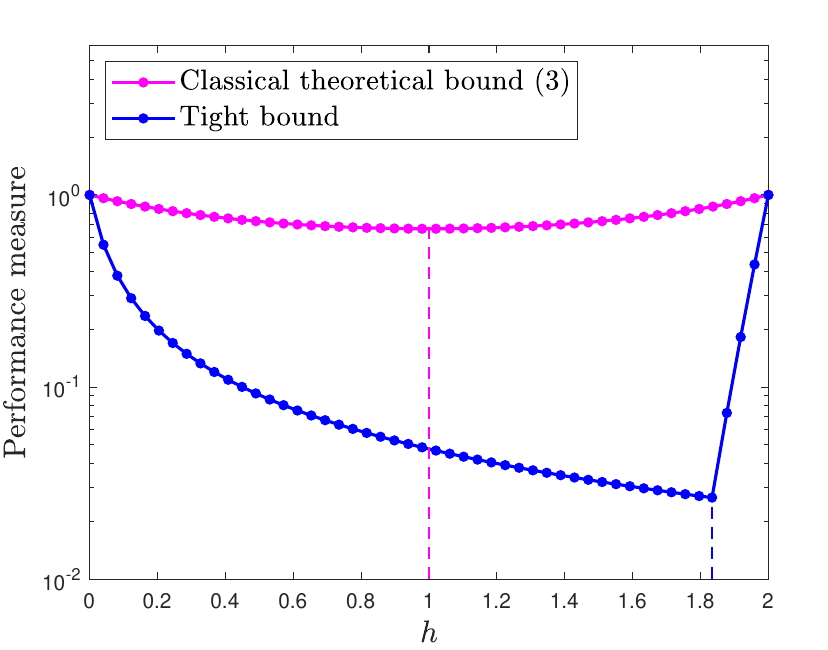}
    \caption{Evolution with the adimensional step-size $h=L\alpha$ of two worst-case guarantees on $f(x_{N})-f(x^*)$ for $L$-smooth convex functions $f$, where the number of iterations $N$ is set to $10$. The plot shows (i) the theoretical bound \eqref{eq:bound_grad_desc} and (ii) the exact bound for this setting. This tight bound allows improving the step-size selection, since the choice $h=1$ minimizing the bound \eqref{eq:bound_grad_desc}  requires doubling the number of steps to obtain the same performance as the true optimal step-size. This highlights the benefits of relying on tight bounds for parameter selection and algorithm comparisons.} \vspace*{-2mm}
\label{fig:compare_grad_desc1_easy}
\end{figure}

A classical way of characterizing the efficiency of optimization is through the derivation of \emph{worst-case bounds}: guarantees on the evolution of a quantity of interest valid for all functions in a given class $\F$.
As a pedagogical example, consider the gradient descent algorithm with a constant step-size $\alpha>0$
\begin{equation}\label{eq:grad_desc}
    x_{i+1}  = x_i - \alpha \nabla f(x_i), \hspace{0.5 cm} \ i=0,\dots,N-1,
\end{equation}
and let $\F_L$ be the set of $L$-smooth convex functions defined on $\R^d$, i.e. functions for which  for all $x,y$ there holds
\begin{align}\label{eq:def_FL}
\begin{cases}
    &f( \lambda x + (1-\lambda)y) \leq \lambda f(x) + (1-\lambda) f(y)   \hspace{.2cm} \forall \lambda\in [0,1], \\&
||\nabla f(x) - \nabla f(y)|| \leq L ||x-y||.
\end{cases}
\end{align}
One can prove, e.g. from \cite[Theorem 2.1.13]{nesterov2003introductory}, that for every function $f\in \F_L$, every $N$ and every $0\leq \alpha\leq \frac{2}{L}$ there holds
\begin{equation}\label{eq:bound_grad_desc}
    f(x_N)-f(x^*) \leq   \frac{2L ||x_0-x^*||^2}{4 + N L\alpha(2-L\alpha) }. \
\end{equation}
Such bounds can be used not only for the intrinsic guarantees they provide, but also to compare the performances of algorithms, and tune their parameters, e.g. $\alpha$ in \eqref{eq:grad_desc}. One would indeed assume that an algorithm with a better bound is more efficient, and typically select the algorithm parameters leading to the best bound. However, this seemingly intuitive approach may lead to sub-optimal or misguided choices when the bounds are too conservative:
an algorithm with a stronger bound is not necessarily more efficient, and may just be more amenable to the analysis technique that was used. Similarly, the parameters minimizing the bound do not necessarily optimize the actual performance, and may just be those for which the proof technique introduces the smallest amount of conservatism. Far from being theoretical possibilities, these phenomena have been observed even in very simple settings.
Fig. \ref{fig:compare_grad_desc1_easy} compares the bound \eqref{eq:bound_grad_desc} with the exact worst-case bound in this setting, obtained via the method described in Section \ref{sec:PEP}.
One can see that the classical bound \eqref{eq:bound_grad_desc} is minimized by $\alpha = L^{-1}$, while the actual performance is minimized by $\alpha \approx 1.834L^{-1}$ \cite{PEP_Smooth}. This step-size improves the bound by a factor 2 compared to $\alpha=L^{-1}$, and can be shown to actually divide the number of iterations required to achieve a given accuracy by 2. Moreover, this phenomenon is accentuated as $N$ increases: the optimal step-sizes approach $2L^{-1}$, a range of values that would seem absurd based on the theoretical bound \cite{PEP_Smooth}.

This demonstrates the importance of bounds that are tight, or as tight as possible, to not only understand the performances of algorithms but also to compare them on a sound base in order to select or tune them, and to guide further progress in algorithmic development.

\subsection{Sources of Conservatism}\label{sec:sources}

At an abstract level, almost all proofs of worst-case bounds consist in (i) translating conceptual assumptions on the function classes into algebraic relations between quantities playing a role in the algorithm, often the iterates $x_i$ and the gradient $\nabla f(x_i)$ and/or function values $f(x_i)$ at these iterates, and (ii) combining these relations with the algorithm description to obtain a bound on a quantity of interest, see e.g. \cite{nesterov2003introductory}.
Conservatism can then come from two possible sources: an insufficiently effective translation of the conceptual assumptions, and/or a sub-optimal combination of the algebraic relations obtained with the algorithm description.

\begin{figure}[h!]
    \centering
    \vspace{-1mm}
    \includegraphics[width=0.48\textwidth]{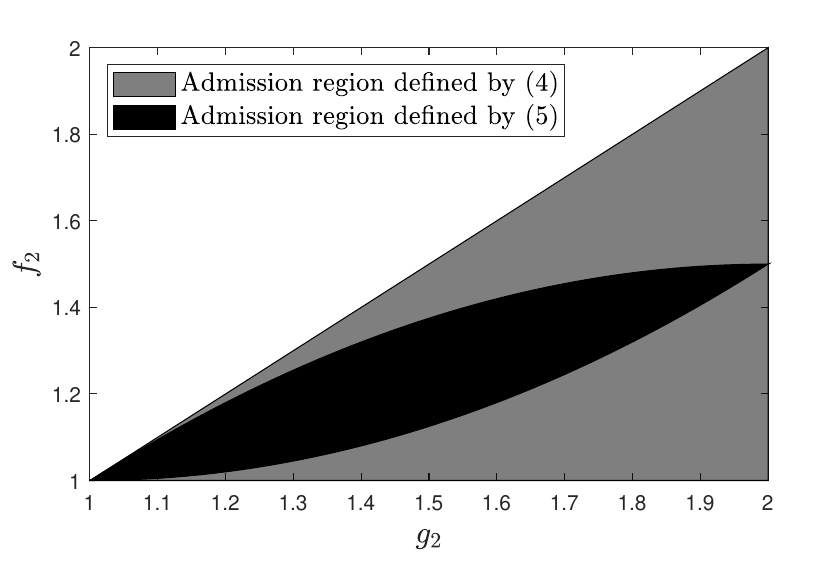}
    \caption{Comparison between the admissible values for $f_2$ and $g_2$ under respectively the discretized versions of \eqref{eq:smoothconvexity} (in grey) and \eqref{eq:interp_smoothconvexity} (in black), assuming $x_2=1,x_1=0,g_1=1,f_1=0$.
    This shows that, despite the equivalence between these constraints when imposed on all pairs $(x,y)$, \eqref{eq:smoothconvexity} is significantly weaker than \eqref{eq:interp_smoothconvexity}. In particular, the points in the grey (but not black) area do not correspond to any actual function $f \in \F_L$.
    \vspace{1cm}}
    \label{fig:nonequivalentdiscretizedconstraints}
\end{figure}

\begin{figure}[h!]
    \vspace*{-3mm}
    \centering
    \includegraphics[width=0.48\textwidth]{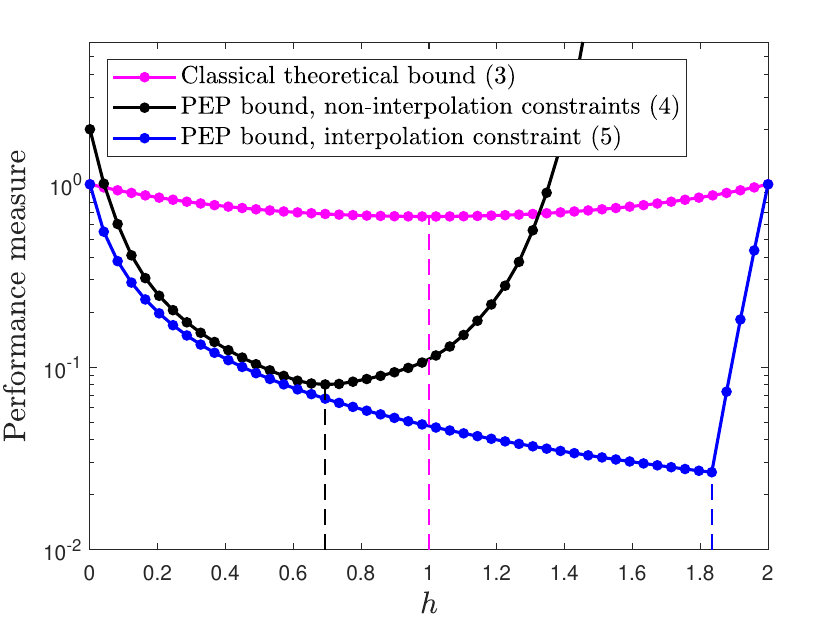}
    \caption{Evolution with the adimensional step-size $h=L\alpha$ of three worst-case guarantees on $f(x_{10})-f(x^*)$ for $L$-smooth convex functions $f$. The plot shows (i) the theoretical bound \eqref{eq:bound_grad_desc}, (ii) the best possible bound (i.e. PEP-based, see Section \ref{sec:PEP}), relying on the non-tight representation of $\mathcal{F}_L$ \eqref{eq:smoothconvexity} (see Section \ref{sec:functions}), and (iii) the tight worst-case bound, i.e. the best possible bound relying on interpolation constraints \eqref{eq:interp_smoothconvexity}. As already observed in Fig. \ref{fig:compare_grad_desc1_easy}, one can see that the step-sizes minimizing non-tight bounds lead to poorer actual performances, even if the bound is PEP-based, that is the best possible bound given a class representation, when the representation is non-tight. This highlights the necessity of relying on tight descriptions of classes to efficiently tune methods.}
    \label{fig:compare_grad_desc1}
\end{figure}

We defer the discussion of the second aspect to Section \ref{sec:PEP} and focus here on the first one, coming back to our example set $\F_L$ of convex $L$-smooth functions
Due to the use of quantifier $\forall \lambda$, \eqref{eq:def_FL} cannot be easily directly manipulated for finite set of points.
However, an equivalent definition of $\mathcal{F}_L$, involving only easily discretizable constraints is the following
\cite{nesterov2003introductory}:
\begin{align}\label{eq:smoothconvexity}
\begin{cases}
    &f(x)\geq f(y)+\langle \nabla f(y),x-y\rangle\\& ||  \nabla f(x)- \nabla f(y)||\leq L||x-y||
\end{cases} \quad  \forall x,y.
\end{align}
An alternative definition, if less intuitive one, is
\begin{equation}\label{eq:interp_smoothconvexity}
    f(x)\geq f(y)+\langle \nabla f(y),x-y\rangle+\frac{1}{2L}||\nabla f(x)-\nabla f(y)||^2,
\end{equation}
$\forall x,y$. One can show that constraints \eqref{eq:smoothconvexity} and \eqref{eq:interp_smoothconvexity} are both completely equivalent to the definition \eqref{eq:def_FL}  of $\F_L$ \emph{when they are imposed for every couple $(x,y)\in \R^d$} \cite{nesterov2003introductory}. However, most proofs will involve a finite number of equations, meaning that the constraints would be imposed on a finite number of pairs $(x,y)$, corresponding typically to the iterates of the algorithm and to remarkable points of the functions. And the \emph{equivalence between the constraints no longer holds when they are imposed on only a finite number of pairs}: Constraint \eqref{eq:smoothconvexity} is indeed weaker than \eqref{eq:interp_smoothconvexity}, and sometimes significantly weaker as demonstrated in Fig. \ref{fig:nonequivalentdiscretizedconstraints}.
Consequently, for most problems, the best bound that can be obtained using \eqref{eq:smoothconvexity} will be weaker than that obtainable from \eqref{eq:interp_smoothconvexity}, as further illustrated in Fig. \ref{fig:compare_grad_desc1}: one can observe that the tightest bound obtained based on constraint \eqref{eq:smoothconvexity} is minimized by a step-size of $\alpha\approx 0.694 L^{-1}$, which yields an even worse actual performance than the classical choice of $\alpha=L^{-1}$, and that the bound explodes as $\alpha$ approaches $2L^{-1}$. 
In other words, relying on \eqref{eq:smoothconvexity} to derive performance guarantees directly introduces more conservatism than using \eqref{eq:interp_smoothconvexity}.

\subsection{Interpolation Constraints}

The simple example in Section \ref{sec:sources} highlights the importance of working with algebraic characterizations of functions classes that introduce the smallest amount of conservatism, or none at all if possible. This is in fact the case of \eqref{eq:interp_smoothconvexity}: it has been shown that for any arbitrary set of $N$ triples $\{(x_i, g_i,f_i)\}_{i\in I}$ satisfying its discretized version
\begin{equation}\label{eq:interp_smoothconvexity_disc}
    f_j\geq f_i+\langle g_i,x_i-x_j\rangle+\frac{1}{2L}||g_i-g_j||^2 \hspace{0.5cm} \forall i,j,
\end{equation}
there exists an actual function $f\in \F_L$ such that $f(x_i)=f_i$ and $\nabla f(x_i) = g_i$ for every $i$. Replacing the assumption $f\in \F_L$ by \eqref{eq:interp_smoothconvexity_disc} being satisfied for all pairs of points appearing in the problem does not introduce any conservatism, since there is always a function in $\F_L$ that would correspond to these values. We say that such a function \emph{interpolates} the $\{(x_i, g_i,f_i)\}$ and that \eqref{eq:interp_smoothconvexity_disc} 
is an \emph{interpolation constraint} since it guarantees the existence of an interpolating function. By contrast, the weaker constraint \eqref{eq:smoothconvexity} is not an interpolation constraint for $\F_L$, as sets of triples $\{(x_i, g_i,f_i)\}$ satisfying it do not necessarily correspond to an actual function $f\in \F_L$, see Fig. \ref{fig:nonequivalentdiscretizedconstraints}. Hence any bound established using \eqref{eq:smoothconvexity} is valid for a class of unspecified \quotes{objects} larger than the functions in $\F_L$ and will thus in most cases be conservative for $\F_L$.

\subsection{Paper Organization}

In this tutorial paper we will first describe how interpolation constraints can be used to automatically compute worst-case bounds in Section \ref{sec:PEP}, and then review interpolation constraints for several classes of objects: functions and general operators in Section \ref{sec:functions}, linear operators in Section \ref{sec:linear}, and doubly-stochastic matrices related to networks in decentralized optimization in Section \ref{sec:network_matrices}.

\section{Automatic Computation of\\Worst-Case Bounds} \label{sec:PEP}
We come back to the second sources of conservatism identified in Section \ref{sec:sources}, i.e. the combination of algebraic translations of the conceptual assumptions and the algorithm description to obtain a bound. Combining such algebraic relations efficiently could in principle be solved by skilled manipulations, but finding the best way of manipulating them can in many cases be very challenging.
Fortunately, this process can be automated, thanks among others to the Performance Estimation Problem (PEP) approach, which was the initial motivation for the formalization of the notion of interpolation constraints in this context \cite{Taylor_thesis,PEP_composite,PEP_Drori,PEP_Smooth} and which we will demonstrate here on our pedagogical example of the gradient descent algorithm \eqref{eq:grad_desc}.

Suppose we want to obtain the best possible bound of the form \eqref{eq:bound_grad_desc}. The idea of PEP is to {compute} the maximal value of $f(x_N)-f(x^*)$ over all possible functions $f\in \F_L$ and starting points $x_0$ for which $||x_0-x^*||^2\leq 1$ (We can then obtain a bound proportional to $||x_0-x^*||^2$ using a simple scaling argument). This a priori infinite-dimensional problem can be re-expressed in a finite manner:
rather than considering the function $f$ as a whole, we will only consider a discrete set $\{\qty(x_i,g_i,f_i)\}_{i\in I= \{0,\dots,N,*\}}$ of iterates and optimal point $x_i$ together with their gradient-vectors $g_i$ and function values $f_i$, \emph{under the constraint that the set is consistent with an actual function $f\in \F_L$}, i.e. that there exists a function $f\in \F_L$ interpolating it.
The optimality of $x^*$ and the fact that the iterates are obtained by the algorithm \eqref{eq:bound_grad_desc} can then also be expressed as constraints. The PEP takes thus the form:
\begin{align}\label{eq:PEP_gd}
  &\underset{\{\qty(x_i,g_i,f_i)\}_{i \in I}}{\max} \quad && f_N - f^* \hspace*{-4mm}  \\[-0.3mm]
 \text{ s.t.}&&  \\
 &\text{algorithm: } &&  x_{i+1} = x_i - \alpha g_i \text{ for $i=0,\dots,N-1.$} \\[-4mm]
 &\text{optimality: } && g^* = 0\\
 &\text{initial cond.: } && ||x_0-x^*||^2\leq 1\\
 &\text{consistency of } && \{\qty(x_i,g_i,f_i)\}_{i \in I} \text{ with some } f\in \F_L.
\end{align}
The optimum of this problem gives by definition a tight worst-case performance bound for the gradient descent with $L$-smooth convex functions. Moreover, the optimal solution corresponds to an instance of function and initial condition actually reaching this upper bound, which can provide very relevant information on the bottlenecks faced by the algorithm. In many cases, these worst-case instances have actually a surprisingly simple structure \cite{Taylor_thesis}.

In \eqref{eq:PEP_gd}, the consistency of $\{\qty(x_i,g_i,f_i)\}_{i \in I} \text{ with some } f\in \F_L$ is implicit, but can be directly replaced by an interpolation constraint when one is available, as \eqref{eq:interp_smoothconvexity_disc} in the case of $\F_L $. When no interpolation constraint is known, one can always relax the consistency constraint into a \emph{necessary constraint for interpolability}, as e.g. \eqref{eq:smoothconvexity}. The solution to the problem obtained will then be a valid worst-case upper bound, and in fact the best upper bound that could possibly be obtained using the necessary constraint used, but in general not a tight one, as illustrated in Fig. \ref{fig:compare_grad_desc1}.

Observe that in both the cases of \eqref{eq:interp_smoothconvexity_disc} and \eqref{eq:smoothconvexity}, the resulting problem is quadratic and potentially non-convex in the iterates and gradients vectors, but it is linear in the scalar products of these and in the function values. Hence, letting the matrix of these scalar products and the vector of function values be the decision variables of the PEP, one can reformulate it as a semi-definite program (SDP) and solve it efficiently, see e.g. \cite[Section 3.2]{PEP_Smooth} for more details. We refer to any constraint linear in these scalar products and function values as \emph{PEP-representable}.

The PEP framework allows computing the exact performance of a wide range of first-order methods on various classes of functions and objects for which PEP-representable interpolation constraints are known. In particular, any method, whether implicit or explicit, consisting in (sub)gradient queries, linear operators, and with PEP-representable constraints can be straightforwardly analyzed.
This includes e.g the gradient method, the proximal gradient method \cite{taylor2018exact}, the coordinate descent \cite{kamri2022worst}, inexact gradient methods \cite{taylor2019stochastic}, Bregman gradient methods \cite{dragomir2021optimal}, variational problems and splitting methods \cite{ryu2020operator,lee2022convergence}, and their combinations. In general, the analysis requires pre-determined step-sizes, possibly time-varying but independent of the function values or gradients. But some efforts have also been made in adaptive first-order methods \cite{de2017worst,barre2020complexity}. The formulation and resolution of Performance Estimation Problems is implemented in Matlab and Python toolboxes:
PESTO \cite{PESTO} and PEPit \cite{pepit2022}. These toolboxes allow for a simple and human readable description of the algorithms and objects, close to their mathematical expression.

Note that interpolation constraints are also used to obtain tight automated convergence guarantees via another approach, relying on tools from robust control \cite{lessard2016analysis, taylor2018lyapunov, upadhyaya2023automated}. The idea is to automatically derive a quadratic Lyapunov function that serves as certificate for the linear convergence of various first-order methods, by analyzing a single step of the method. Provided that the method is analysed over a class of functions for which interpolation constraints are available, the Lyapunov function obtained is tight in the sense that it guarantees the best decrease rate possible for any \emph{single} (or fixed number of) iteration of the algorithm and for a quadratic Lyapunov function. Compared to the PEP framework, one of the assets of this approach is its tractability, since the SDP to solve is of small size. One the other hand, it cannot deal with time-varying steps or sublinear rates, and the PEP framework can achieve better rates by analyzing several steps.

\section{Functions and Operators Interpolation} \label{sec:functions}
In this section, we review the functions classes $\mathcal{F}$ for which the following question has been answered: what necessary and sufficient constraints, i.e. interpolation constraints, must a set of data (points, functions values and (sub)gradients) satisfy to ensure its $\mathcal{F}$-interpolability:
\begin{definition}[$\mathcal{F}$-interpolability]\label{def:func_interp}
    Let $I$ be a finite set of indices. A set of triples $\{(x_i,g_i,f_i)\}_{i\in I}$ is $\mathcal{F}$-interpolable if
    \begin{equation}
        \exists F\in \mathcal{F}~:~
        \begin{cases}
            f_i = F(x_i) \\
            g_i = \nabla F(x_i)
        \end{cases}
         \forall i\in I.
    \end{equation}
\end{definition}

The concept of interpolation constraint applied to optimization was introduced in \cite{PEP_Smooth}, to tightly analyze first-order methods on the class of smooth strongly convex functions $\mathcal{F}_{\mu,L}$, classically defined by:
\begin{align}\label{eq:smoothstrongconvexity}
\hspace{-5mm}
\begin{cases}
    f(x)\geq f(y)+\langle \nabla f(y),x-y\rangle+\frac{\mu}{2}||x-y||^2\\
    ||  \nabla f(x)- \nabla f(y)||\leq L||x-y||
\end{cases} \forall x,y.
\end{align}
When $\mu=0$, definition \eqref{eq:smoothstrongconvexity} amounts to definition \eqref{eq:smoothconvexity} of $\mathcal{F}_L$. Separately, the discretized versions of the two constraints in \eqref{eq:smoothstrongconvexity} are interpolation constraints: for instance, a data set  $\{(x_i,g_i,f_i)\}_{i\in I}$ is interpolable by a strongly convex function if and only if it satisfies
\begin{align}\label{eq:strongconvexity}
    &f_i\geq f_j+\langle g_j,x_i-x_j\rangle+\frac{\mu}{2}||x_i-x_j||^2 \quad \forall i,j\in I.
\end{align}
However, it can be shown that the juxtaposition of these two interpolation constraints, i.e. the discretized version of Definition \eqref{eq:smoothstrongconvexity}, is not tight, while an interpolation constraint for $\mathcal{F}_{\mu,L}$ is given by :
\begin{theorem}[{$\mathcal{F}_{\mu,L}$}-interpolation constraint \cite{PEP_Smooth, Rubbens_interp}]\label{th:int_cond_mul}
 A set of triples $\{(x_i,g_i,f_i)\}_{i\in I}$ is $\mathcal{F}_{\mu,L}$-interpolable if
\begin{align} \label{cond:cond_mul}
f_i&\geq f_j+\langle g_j,x_i-x_j\rangle +\frac{||g_i-g_j||^2}{2(L-\mu)}\\
&+\frac{\mu L||x_i-x_j||^2}{2(L-\mu)}-\frac{\mu}{L-\mu}\langle g_i-g_j,x_i-x_j\rangle \ \forall i,j,
\end{align}
when $\mu\neq L$, and
\begin{align} \label{cond:cond_mu2l}
f_i&\geq f_j+\frac{1}{2}\langle g_i+g_j,x_i-x_j\rangle \\&+\frac{||g_i-g_j-L(x_i-x_j)||^2}{L}\  \forall i,j,
\end{align}
otherwise.
\end{theorem}
Hence, an analysis relying on these interpolation constraints is a priori tighter than any analysis relying on different definitions of this function class, as shown in Fig. \ref{fig:compare_grad_desc1} for $\mu=0$. It was later shown in \cite{rotaru2022tight} that this result also holds for negative values of $\mu$ as well, that is for smooth \emph{weakly convex} functions, i.e. convex up to the addition of a quadratic term. Interpolation constraints for $\mathcal{F}_{\mu,L}$ that do not explicitly imply function values, in case one would not want to use these as variables, also exist, based on Rockafellar's cyclic monotonicity conditions \cite{rockafellar1997convex, Taylor_thesis}:
\begin{theorem}[{$\mathcal{F}_{\mu,L}$}-interpolation without function values]\label{th:int_cond_mul_nof}
 A set of pairs $\{(x_i,g_i)\}_{i\in I}$ is $\mathcal{F}_{\mu,L}$-interpolable, in the sense that $\exists F\in \mathcal{F}_{\mu,L}~:~ g_i = \nabla F(x_i)\quad    \forall i\in I$,
  if for any cyclic sequence $(x_{i_1},g_{i_1}), \dots, (x_{i_N},g_{i_N}), (x_{i_1},g_{i_1})$,
\begin{align}
\sum_{k=1}^N &\langle g_{i_k},x_{i_k}-x_{i_{k+1}}\rangle -\frac{\langle g_{i_k},g_{i_k}-g_{i_{k+1}}\rangle }{L}\\
&-\mu \langle x_{i_k},x_{i_k}-x_{i_{k+1}}\rangle-\frac{\mu \langle x_{i_k},g_{i_k}-g_{i_{k+1}}\rangle }{L}\geq 0. \\[-4mm]
\end{align}
\end{theorem}
However, checking these constraints rapidly becomes intractable when the data set size increases, so that when possible, introducing function values as variables is more efficient.

In \cite{PEP_Smooth}, interpolation constraints are also given for smooth functions with bounded subgradient $\mathcal{C}_{L,M}$ and strongly convex functions with bounded domain $\mathcal{S}_{\mu,D}$. While, in general, the juxtaposition of two interpolation constraints does not remain an interpolation constraint, it is the case in these particular cases:
\begin{theorem}[{$\mathcal{C}_{L,M}$}-interpolation constraints \cite{Taylor_thesis}]\label{th:int_cond_clm}
 A set of triples $\{(x_i,g_i,f_i)\}_{i\in I}$ is $\mathcal{C}_{L,M}$-interpolable if
\begin{align} \begin{cases}
&f_i\leq f_j+\langle g_j,x_i-x_j\rangle +\frac{L}{2}||x_i-x_j||^2\\
&||g_i||\leq M.
\end{cases}
\end{align}
\end{theorem}

On the other hand, one can show that this property does not hold anymore for the class of weakly convex functions with bounded subgradient, for which the juxtaposition of the two constraints is not interpolable \cite{Rubbens_interp}. Finally, interpolation constraints are also available for indicator functions \cite{Taylor_thesis}:
\begin{align}
    i_C(x)=\begin{cases}
        &0 \text{ if } x\in C\\
        & 1 \text{ otherwise}.
    \end{cases}
\end{align} for $C$ a closed convex set. Let $\mathcal{I}_M$ be the class of indicator functions bounded in radius by some M, i.e. $||x||\leq M \ \forall x \in C$. These functions are of especially high importance for projections in constrained optimization, and their interpolation constraints are given by:
\begin{theorem}[{$\mathcal{I}_M$}-interpolation constraints \cite{Taylor_thesis}]\label{th:int_cond_indic}
 A set of triples $\{(x_i,g_i,f_i)\}_{i\in I}$ is $\mathcal{I}_M$-interpolable if
\begin{align} \begin{cases}
&f_i=0\\&
\langle g_j,x_i-x_j\rangle\leq 0\\
&||x_i||\leq M.
\end{cases}
\end{align}
\end{theorem}
\smallskip
Aside from these interpolation constraints for function classes, we can adapt Definition \ref{def:func_interp} to consider operator classes $\mathcal{Q}$. In this work, we focus on real operators which are defined as multi-dimensional mappings $Q:\mathbb{R}^d\rightarrow\mathbb{R}^m$.
\begin{definition}[$\mathcal{Q}$-interpolability]\label{def:op_interp}
    Let $I$ be a finite set of indices. A set of pairs $\{(x_i,q_i)\}_{i\in I}$ is $\mathcal{Q}$-interpolable if
    \begin{equation}
        \exists Q\in \mathcal{Q}~:~
            q_i = Q(x_i) \quad
         \forall i\in I.
    \end{equation}
\end{definition}
Note that operator interpolation can be viewed as zero-th order function interpolation for multi-valued functions. Moreover, the concept of interpolation without function values, e.g Theorem \ref{th:int_cond_mul_nof}, can be seen as operator interpolation for a subclass of operators, coherent as gradients of a function. Interpolation constraints exist for (strongly) monotone $\mathcal{M}_{\mu}$, cocoercive $\mathcal{C}_{\beta}$ and Lipschitz operators $\mathcal{N}_L$\cite{whitney1992analytic, bauschke2011convex, ryu2020operator}:
\begin{theorem}[{ $\mathcal{M}_{\mu}$,  $\mathcal{C}_{\beta}$, $\mathcal{N}_{L}$-interpolation \cite{ryu2020operator}}]\label{th:int_op}
 A set of pairs $\{(x_i,q_i)\}_{i\in I}$ is
\begin{align}
&\mathcal{M}_{\mu}\text{-interpolable if }\forall i,j, \ \langle q_i-q_j,x_i-x_j\rangle\geq \mu ||x_i-x_j||^2,\\
&\mathcal{C}_{\beta}\text{-interpolable}\quad \text{if }    \forall i,j, \  \langle q_i-q_j,x_i-x_j\rangle\geq \beta ||q_i-q_j||^2,\\
&\mathcal{N}_{L}\text{-interpolable }\  \text{if }\forall i,j, \ ||q_i-q_j||\leq L ||x_i-x_j|| .
\end{align}
\end{theorem}
However, no interpolation result is known for operators combining several of these properties.

As explained in Section \ref{sec:PEP}, combining the PEP framework with any of these interpolation constraints, when PEP-representable (i.e. linear in the function values and scalar products of points and (sub)gradients), allows tightly analyzing a wide variety of first-order methods on these function and operator classes. Even when interpolation constraints are not PEP-representable, involving e.g. absolute values or exponentials, it holds that any theoretical analysis based on them will be a priori tighter than an analysis based on non interpolable constraints, which explains why deriving interpolation constraints has an interest far beyond the PEP framework.

\section{Linear Operator Interpolation}\label{sec:linear}

\subsection{Motivations}
Many optimization problems involve a particular subclass of operators that has not yet been described in Section \ref{sec:functions}, that is linear operators $\mathcal{Q}(x) = Mx$. These problems, including a large number of classical optimization problems as total variation deblurring, basis pursuit, resource allocation, etc., contain the linear operator in their objective function, e.g.,
\begin{equation}
    \min_x f(x) + h(Mx),
\end{equation}
or in their constraints
\begin{equation}
    \min_x f(x) \text{ s.t. } Mx = b,
\end{equation}
possibly with more terms in the objective function or constraint (e.g. $f(x) + g(M_1x) + h(M_2x)$ or $M_1 x + M_2y = b$).

Consequently, many methods have been designed to solve such structured problems, e.g. Chambolle-Pock method \cite{chambolle2011first}, primal-dual fixed point method \cite{chen2016primal} or alternating direction method of multipliers \cite{gabay1976dual}.
They typically combine gradient and linear operator queries. Therefore, tightly analyzing these methods in PEP requires obtaining interpolation constraints representing a linear operator or, more precisely, representing its effect.

In this section, we consider finite-dimensional linear operators, i.e. matrices. The matrices considered are not necessarily symmetric and have their maximal singular value bounded by some $L\geq 0$. In other words, we work with the following set of linear operators:
\begin{align}
\begin{split}
    \LO{L} & = \{ M: \sigma_{\max}(M) \leq L\}.
\end{split}
\end{align}

\subsection{Interpolation Constraints for Linear Operators}
To allow tightly analyzing methods involving linear operators, we now explain how to represent the application of a linear operator in PEP.

Let us consider a function $F(x) = h(Mx)$, where $h\in \F$, and its (sub)-gradient $\nabla F(x_i) = M^T \nabla h(Mx_i)$ at $x_i$. Introducing intermediate variables $y_i$, $u_i$ and $v_i$ allows decomposing the gradient as
\begin{align}
    \begin{split}
        y_i & = Mx_i, \\
        u_i & = \nabla h(y_i), \\
        v_i & = M^T u_i = \nabla F(x_i).
    \end{split}
\end{align}

If interpolation constraints are known for $\F$, the existence of some $h\in \F$ such that $u_i = \nabla h(y_i)$ is guaranteed by imposing the corresponding interpolation constraints on the set $\{(y_i,u_i)\}_{i \in I}$, involving potentially unused function values $h_i$ if necessary, see Section \ref{sec:functions} for more details.
To impose constraints $y_i = Mx_i$ and $v_i = M^Tu_i$, one could consider $M$ as a variable of the PEP and try to impose conditions on its singular values. However, the problem would then become non convex since involving the multiplication of variables. We rather take the view of interpolation constraints for operators
 and impose convex constraints on the sets $\{(x_i,y_i)\}$ and $\{(u_i,v_i)\}$, guaranteeing the existence of a matrix $M\in \mathcal{L}_L$ interpolating them.

\begin{definition}[$\mathcal{L}_L$-interpolability]\label{def:R_mat_inter1}
 Let $I=\{1,\ldots,N\}$ a set of indices. Sets of pairs $\{(x_i,y_i)\}_{i\in I}$ and $\{(u_i,v_i)\}_{i\in I}$ are $\LO{L}$-interpolable if
    \begin{equation}
        \exists M\in \LO{L}~:~
        \begin{cases}
            y_i = M x_i, & \forall i\in I, \\
            v_i = M^T u_i, & \forall i\in I.
        \end{cases}
    \end{equation}
\end{definition}
In \cite{bousselmi2023interpolation}, interpolation constraints for the class $\mathcal{L}_L$ are derived:
\begin{theorem}[{$\LO{L}$}-interpolation constraints]\label{th:int_cond_non_sym}
 Sets of pairs $\{(x_i,y_i)\}_{i\in I}$ and $\{(u_i,v_i)\}_{i\in I}$ are $\LO{L}$-interpolable if
\begin{align} \label{cond:non_symm}
\begin{split}
\begin{cases}
        X^T V = Y^TU, \\
        Y^TY \preceq L^2 X^TX,\\
        V^TV \preceq L^2 U^TU,
\end{cases}
\end{split}
\end{align}
where $X=(x_1~\dots~x_{N})$, $Y=(y_1~\dots~y_{N})$, $U=(u_1~\dots~u_{N})$ and $V=(v_1~\dots~v_{N})$.
\end{theorem}

Observe that these constraints are convex semidefinite constraints on the scalar products of the points and therefore PEP-representable.
They can hence be used to analyze worst-case performance of any method involving linear operators via PEP. In \cite{bousselmi2023interpolation}, the authors exploit this result to analyze, e.g., the Chambolle-Pock method \cite{chambolle2011first} or the gradient method on a structured function $h(Mx)$, where $h\in F$ and $\F$ is a function class for which interpolation constraints are known. Note that Theorem \ref{th:int_cond_non_sym} can be extended to symmetric linear operators with bounded eigenvalues, see \cite{bousselmi2023interpolation}.

\section{Network Matrices Interpolation}\label{sec:network_matrices}

In this section, we focus on the interpolation of a specific set of finite-dimensional linear operators: the set of symmetric doubly-stochastic matrices, which have a given range of eigenvalues.
These matrices are often used to represent consensus steps over a network of agents, which are a building block of most distributed optimization methods. Therefore, having such interpolation constraints allows building a PEP computing worst-case guarantees for distributed optimization algorithms.
The main difference with the set of linear operator $\LO{L}$ characterized in Section \ref{sec:linear} is the stochasticity of the matrix, required for having an averaging consensus.

\subsection{Distributed Optimization}
Distributed optimization methods aim at minimizing the average of local functions that are distributed across a network of agents $\{1,\dots,\A\}$. They are used, for example, in large-scale machine learning or sensor networks for various reasons, including data privacy, easy scaling to large size problems, robustness to failure and possible speed up.
Formally, the goal of decentralized optimization algorithms is to collaboratively solve the following optimization problem:
\begin{equation*}
 \underset{\text{\normalsize $x \in\mathbb{R}^{d}$}}{\mathrm{minimize}} \quad f(x) = \frac{1}{\A}\sum_{\a = 1}^{\A} f_{\a}(x),  \vspace{-1mm}
\end{equation*}
where $f_\a: \mathbb{R}^{d}\to\mathbb{R}$ is the private function locally held by agent $\a$.
Each agent $\a$ holds its own version $x_\a \in \Rvec{d}$ of the decision variable $x$, performs local computations and exchanges local information with its neighbors to come to an agreement on the minimizer $x^*$ of the global function $f$.
In this work, we focus on decentralized algorithms where the exchanges of information take the form of an average consensus on some quantity, e.g.,  on the $x_\a$. One of the first and simplest algorithm of this form is the \emph{distributed gradient descent} (DGD) \cite{DGD1,DGD} for which iteration $i$ can be described as \vspace{-1mm}
\begin{align}
  y_{a,i} &= \sum_{b=1}^\A w_{ab} \, x_{b,i},          & \text{ for $\a = 1,\dots, \A.$} \hspace*{5mm} \label{eq:DGD_cons} \\[-0.5mm]
  x_{\a,i+1} &= y_{\a,i} - \alpha_i \nabla f_{\a}(x_{\a,i}), & \text{ for $\a = 1,\dots, \A.$} \hspace*{5mm} \label{eq:DGD_comp} \vspace{-2mm}
\end{align}
The vectors $y_{\a,i} \in \Rvec{d}$ represent the result of the consensus step. The matrix $W \in \Rmat{\A}{\A}$ contains all the averaging weights $w_{ab}$ and is called the network or averaging matrix.
Consensus step \eqref{eq:DGD_cons} can be written as a matrix multiplication where the variables of all the agents are stacked in single columns $\x_i, \y_i \in \Rvec{\A d}$:
\begin{equation}
    \y_i = (W \kron I_d) \x_i, \quad \text{with} \quad \y_i = \begin{bmatrix} y_{1,i} \\ \vdots \\ y_{\A,i} \end{bmatrix}   \text{ and }  \x_i = \begin{bmatrix} x_{1,i} \\ \vdots \\ x_{\A,i} \end{bmatrix},
\end{equation}
where $\kron$ denotes the Kronecker product and $I_d$ the identity matrix of size $d$.
This corresponds to applying a linear operator with a particular matrix shape and properties.
In what follows, we consider the following set of network matrices, for a fixed number of agents $\A > 1$ and a fixed bound $\lam \in [0,1)$ on the non-principal eigenvalues:
\begin{equation}
    \Wcl{\lam} = \left\{ W \in \Rmat{\A}{\A}: \parbox[c]{5.2cm}{ \text{$W$ is symmetric and $\lam_1(W) = 1$,} \\ \text{$-\lam \le \lam_\A(W) \le \dots \le \lam_2(W) \le \lam$}} \right\}.
\end{equation}
The performance guarantees from the literature, see \cite{DGD} for a survey, usually add a non-negativity assumption for the network matrices, by qualifying them as \emph{doubly-stochastic}, but these guarantees are in fact often valid for all the matrices from the set $\Wcl{\lam}$.
Representing $\Wcl{\lam}$ in PEP through interpolation constraints is therefore important to obtain general and comparable bounds for DGD or any other decentralized algorithm involving consensus steps.

\subsection{Interpolation Constraints for Consensus Steps}

We now present interpolation constraints to describe the effect of consensus steps with an averaging matrix from $\Wcl{\lam}$. This allows to formulate the performance estimation problem for any decentralized algorithm using averaging consensuses. To analyze such algorithms with the PEP framework, we need to find the worst-case for the local function and the sequence of local iterates of each agent. The same techniques as in Section \ref{sec:PEP} can be applied to discretize the problem, using proper interpolation constraints on each local function. In this section, we focus on the representation of the consensus steps that are part of these algorithms: \vspace{-1mm}\begin{equation}
    \y_i = (W \kron I_d) \x_i \quad \text{ with $W \in \Wcl{\lam}$.}
\end{equation}
As in Section \ref{sec:linear}, we cannot add matrix $W$ as variable of the PEP problem because it makes it intractable but we rather impose convex constraints on the (scalar products of the) iterates $\y_i$ and $\x_i$ that are necessary for the existence of the desired network matrix $W$. The idea is thus to obtain a relaxed representation of the network matrix by formulating necessary constraints on the discretized points. As explained in Section \ref{sec:PEP}, sufficiency of the formulation can be verified a posteriori and will be discussed below.

\begin{definition}[$\Wcl{\lam}$-interpolability] Let $I$ be a set of indices of consensus steps. A set of pairs $\{(\x_i,\y_i)\}_{i\in I}$ is $\Wcl{\lam}$-interpolable if,
    \begin{equation}
        \exists W\in \Wcl{\lam}~:~ \y_i = (W \kron I_d) \x_i \quad \text{ for $i \in I.$}
    \end{equation}
\end{definition}
\smallskip
While we do not have necessary and sufficient interpolability constraints for $\Wcl{\lam}$, Theorem \ref{thm:int_cond_consmat} below provides necessary conditions that have been shown to lead to good results in practice. It is expressed in terms of the matrices ${X=(\x_1~\dots~\x_{N})}$, $Y=(\y_1~\dots~\y_{N})$, as well as their decomposition in two terms, the average and centered parts:
$$ X = (\Xb \otimes \mathbf{1}_\A ) + \Xp, \qquad Y =  (\Yb \otimes \mathbf{1}_\A) + \Yp,$$
where $\Xb$ and $\Yb$ are agents average vectors in $\Rmat{d}{N}$, defined as $\Xb_{\cdot i} = \frac{1}{\A} \sum_{\a=1}^\A x_{\a,i}$, $\Yb_{\cdot i} = \frac{1}{\A} \sum_{\a=1}^\A y_{\a,i}$ for $i=1,\dots,N$ and $\mathbf{1}_\A = \qty[1\dots 1]^T \in \Rvec{\A}$. 

\begin{theorem}[{$\Wcl{\lam}$}-interpolation necessary constraints \cite{PEP_dec}]\label{thm:int_cond_consmat}
 If a set of pairs $\{(\x_i,\y_i)\}_{i\in I}$ is $\Wcl{\lam}$-interpolable, then it holds that
    \begin{align}
    &\text{(Average preservation)} \quad & \Xb &= \Yb, \hspace{1cm}\label{eq:eq_avg} \\
    &\text{(Symmetry)} \quad & \Xp^T\Yp &= \Yp^T\Xp, \hspace{1cm}\label{eq:sym_xyc0} \\
    &\text{(Variance reduction)} \quad & \Yp^T\Yp &\preceq \lam^2 \Xp^T\Xp. \hspace{1cm} \label{eq:var_yyc0}
 \end{align}
\end{theorem}
\smallskip
Constraints \eqref{eq:eq_avg} of average preserving come from the double-stochasticity of the network matrix $W$; the symmetry constraint \eqref{eq:sym_xyc0} directly reflects the symmetry of $W$; and the constraint \eqref{eq:var_yyc0} of variance reducing is a consequence of the bound $\lam$ on the eigenvalues of $W \in \Wcl{\lam}$. All these constraints can be expressed as linear constraints in terms of scalar products of the iterates and are therefore PEP-representable, allowing to automatically compute performance guarantees of a given distributed optimization methods. The guarantees obtained are valid for any averaging matrix $W \in \Wcl{\lam}$. This theorem can be extended to general range of non-principal eigenvalues $[\lam^-, \lam^+] \ne [-\lam, \lam]$ with $\lam^- \le \lam^+ \in (-1,1)$, see \cite[Theorem 1]{PEP_dec}.

In contrast with previous sections, the constraints from Theorem \ref{thm:int_cond_consmat} are necessary, but not sufficient to the $\Wcl{\lam}$-interpolability, which remains an open question. This means that these constraints describe a feasible set containing at least all the sets of iterates $X, Y$ that are $\Wcl{\lam}$-interpolable, but maybe also some others, for which there exists no matrix $W \in \Wcl{\lam}$ such that $Y = (W \kron I_d) X$. Therefore, using these constraints in a PEP provides a valid upper bound on the worst-case performance of the algorithm we analyze. If the worst-case solution obtained with PEP corresponds to iterates for which there is a matrix $W^* \in \Wcl{\lam}$ such that $Y = (W^* \kron I_d) X$, then we know that the worst-case bound is exact.
Moreover, based on the results of Section \ref{sec:linear}, one can show that we can always recover a matrix $M \in \Rmat{\A d}{\A d}$ which is symmetric and with $\lam_1(M) = 1$, $-\lam \le \lam_\A(M) \le \ldots \le \lam_2(M) \le \lam$ but for which we cannot guarantee the shape $M = W \kron I_d$.

\begin{figure}[h!]
  \centering
  \includegraphics[width=0.5\textwidth]{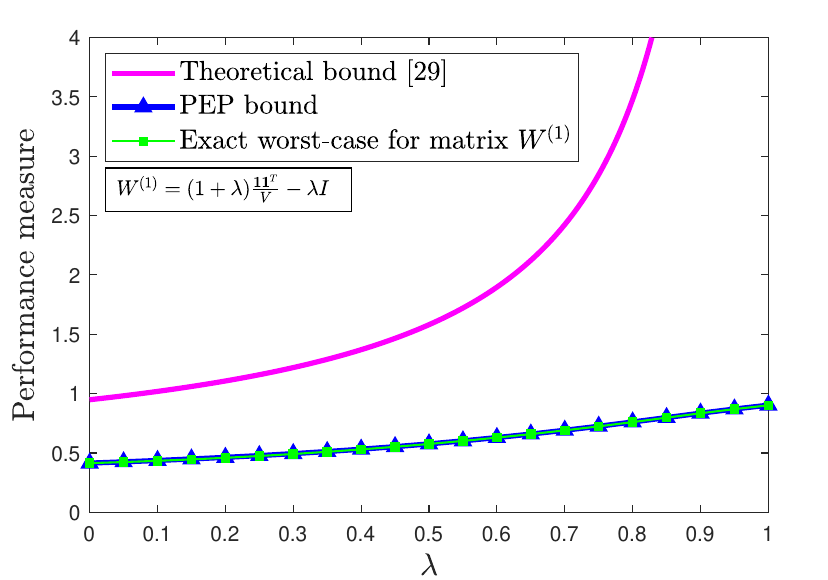}
  \vspace{-3mm}
  \caption{This plot is is inspired from \cite[Fig. 2]{PEP_dec} and shows the evolution with $\lam$ of the worst-case performance of $N = 10$ iterations of DGD with $V = 3$ agents and ${\alpha = \frac{1}{\sqrt{N}}}$. The plot shows (i) the theoretical bound from \cite{DGD} (in pink), largely above (ii) the worst-case performance obtained with PEP using constraints from Theorem \ref{thm:int_cond_consmat}  (in blue) and (iii) the exact worst-case performance for the averaging matrix $W^{(1)}$  (in green).
  The blue and green curves are matching, which indicates the tightness of the spectral PEP bound for DGD.
  \vspace{-2mm}}
  \label{fig:wc_lamevol_N3}
\end{figure}

Another way of verifying the exactness of the PEP bound is to obtain the same worst-case value using an exact PEP problem in which the network matrix $W$ is fixed a priori, with some guess of worst-case network matrix $W^*$. This technique allowed us to verified exactness of the PEP bounds for different algorithms in our different experiments \cite{PEP_dec}, \cite{PEP_dec_CDC} and \cite{PEP_dec_CDC22}. This is shown, for example, for 10 iterations of DGD in Fig. \ref{fig:wc_lamevol_N3}. The figure also shows the large improvement achieved by the PEP bound with respect to the theoretical bound from the literature for constant step-size \cite{DGD}, especially for large values of $\lam$, corresponding to situations where the network is poorly connected.

\section{Conclusions}
We have shown how choices made when translating conceptual assumptions on functions, operators etc. into algebraic relations between quantities involved in an algorithm could introduce significant conservatism in the derivation of worst-case bounds on optimization algorithm, before even beginning to combine these relations. This motivates the study of interpolation constraints: algebraic relations that do not introduce any degree of conservatism because they are necessary and sufficient for the existence of a function or operator consistent with the quantities analyzed.

We have reviewed the interpolation constraints for different classes of functions and general operators, but also for a spectral class of linear operators, that are used in many methods designed to solve structured problems, and for a spectral class of symmetric and stochastic matrices, that are typically used in distributed optimization algorithms.

By definition, interpolation constraints allow in principle obtaining tight bounds on the performance of any optimization algorithm provided one can combine them in an appropriate, or optimal, manner. This is the goal of performance estimation problems (PEP), which obtain bounds by solving an optimization problem of which the algebraic descriptions of the objects at stakes and the algorithmic description are constraints, and which lead to (numerical) tight bounds if interpolation constraints are used. PEPs motivated the study of interpolation constraints, but we stress that their benefit are also valid in classical analysis, for instance when analyzing function classes or methods that are not PEP representable, and advise therefore always using interpolation constraints when they are available.

There remain, however, several challenges. First, interpolation constraints are only available for a subset of the interesting classes of functions and operators. To the best of our knowledge, no interpolation constraint is known yet for e.g. blockwise smooth functions \cite{kamri2022worst}, Holder smooth functions \cite{bonic1966smooth}, smooth functions satisfying a Polyak-Łojasiewicz (PL) condition \cite{karimi2016linear}, or even functions that are smooth and convex but only defined on a subset of $\R^d$, and for which condition \eqref{eq:interp_smoothconvexity_disc} has been shown to be too strict \cite{drori2018properties}. Available interpolation constraints also concern mostly first-order descriptions;  there are no tractable interpolation constraints combining information about first and second-order derivatives, preventing their use in the analysis of second-order methods.
Second, we have no principled way of combining interpolation constraints to tightly describe function classes defined by several properties. Indeed, it is known that simply juxtaposing the interpolation constraints for two classes of functions, e.g. smooth and convex functions, does in general not provide an interpolation constraint for the intersection of these classes, e.g. smooth convex functions. Finally, understanding the actual impact of using approximate (or only necessary) interpolation constraints also remains an open question.

\bibliographystyle{IEEEtran}
\bibliography{refs.bib}
\end{document}